\documentclass{article}      

\usepackage{epsfig}

\title{q-Bernoulli Numbers and Zeros of $q$-Sine Function}

\author{Sengul Nalci and Oktay K. Pashaev \\Department of Mathematics, Izmir Institute of Technology \\ Urla-Izmir, 35430, Turkey}

\begin{document}
\newcommand{\be}{\begin{equation}}
\newcommand{\ee}{\end{equation}}
\newcommand{\bea}{\begin{eqnarray}}
\newcommand{\eea}{\end{eqnarray}}
\newcommand{\disp}{\displaystyle}
\newcommand{\la}{\langle}
\newcommand{\ra}{\rangle}

\newtheorem{thm}{Theorem}[subsection]
\newtheorem{cor}[thm]{Corollary}
\newtheorem{lem}[thm]{Lemma}
\newtheorem{prop}[thm]{Proposition}
\newtheorem{defn}[thm]{Definition}
\newtheorem{rem}[thm]{Remark}
\newtheorem{prf}[thm]{Proof}

\maketitle


\begin{abstract}
There exists a well-known relation between the zeros of sine function, Bernoulli numbers and the Riemann Zeta function. In the present paper, we find a similar relation for zeros of q-sine function. We introduce a new q-extension of the Bernoulli numbers with generating function written in terms of both Jackson's q-exponential functions. By q-generalized multiple product Leibnitz rule and the q-analogue of logarithmic derivative we established exact relations between zeros of $\sin_q x$ and our q-Bernoulli numbers. These relations could be useful for analyzing approximate and asymptotic formulas for the zeros and solving BVP for $q$-Sturm-Liouville problems.
\end{abstract}

\section{Introduction}
One of the most impressive applications of Bernoulli numbers is related with zeros of $\sin x$
function and the Riemann Zeta function. In the present paper, by proper generalization of Bernoulli numbers to $q$-Bernoulli numbers, we establish similar relation for zeros of $\sin_q x$ function. The generating functions of our Bernoulli polynomials and Bernoulli numbers are defined in terms of Jackson's q-exponential functions. Our generating functions and Bernoulli numbers are different from the known in literature \cite{Al-Salam}, \cite{Carlitz}, \cite{Ryoo}. By q-differentiation of the generating function we get the recursion formula for q-Bernoulli polynomials which is reducible to the standard recursion formula in the limit $q \rightarrow 1.$ From the definition of q-exponential functions we write the power series expansion of this generating function and obtain first few q-Bernoulli numbers. In order to write
$\sin_q x$ function as an infinite product in terms of its zeros, we introduce the q-generalized multiple product Leibnitz rule and the q-analogue of logarithmic derivative. This gives us relation between zeros of q-sine function and q-Bernoulli numbers. In the limit  $q\rightarrow 1,$ our results transforms to the known relations between the zeros of sine function, Bernoulli numbers and Riemann Zeta function.

\section{Zeros of Sine Function and Riemann Zeta Function}
First we briefly review the known relation between the zeros of $\sin x$ function, Bernoulli numbers and the Riemann Zeta function.
\subsection{Bernoulli Polynomials and Numbers}
The generating function
\be F_x(z)= \frac{z e^{z x}}{e^z-1} = \sum_{n=0}^\infty B_n(x) \frac{z^n}{n!} \label{bernoullipoly} \ee
 determines the Bernoulli polynomials in $x$,\,$B_n(x),$ $\forall n>0.$

By differentiating in $x$ we get the recursion formula for Bernoulli polynomials

\bea B_n'(x)=n B_{n-1}(x),\,\,\,\,\,\,\, n \geq 1. \label{differentialrelation}\eea
In addition, we have
\bea \forall n \geq 1, \,\,\, B_n (x+1)-B_n(x)= n x^{n-1}. \eea

\textbf{Bernoulli numbers} are defined as $B_n(0)=b_n.$
Then the generating function for Bernoulli numbers follows from (\ref{bernoullipoly})
\be \frac{z}{e^z-1}= \sum_{n=0}^\infty b_n \frac{z^n}{n!}.\label{bernoullinumber}\ee
Below we display first few Bernoulli polynomials and numbers
$$B_0(x)=1,\,\,\,\,\, B_1(x)= x-\frac{1}{2},\,\,\,\,\, B_2(x)=x^2-x+\frac{1}{6},\,\,\,\,\,\, B_3(x)= x^3-\frac{3}{2}x^2+\frac{1}{2}x.$$
$$b_0=1,\,\,\,\,\, b_1=-\frac{1}{2},\,\,\,\,\, b_2= \frac{1}{6},\,\,\,\, b_3=0.$$

Bernoulli numbers are related with zeros of $\sin z$ function, and allows one to calculate the values of the Riemann Zeta function at even numbers argument \cite{Sury}. We consider infinite product representation for $\sin z:$
\be \sin z= z \prod_{n=1}^\infty \left(1-\frac{z^2}{\pi^2 n^2}\right).\ee
\bea  \frac{d}{d z} \ln(\sin z)&=& \frac{d}{d z}\left(\ln \left(z\prod_{n=1}^\infty \left(1-\frac{z^2}{\pi^2 n^2}\right) \right)\right)= \frac{d}{d z} \left(\ln z+ \ln \sum_{n=1}^\infty \left(1-\frac{z^2}{\pi^2 n^2}\right)\right) \nonumber \\
&=& \frac{\cos z}{\sin z}= \frac{1}{z}+ \sum_{n=1}^\infty \frac{\frac{-2 z}{\pi^2 n^2}}{1-\frac{z^2}{\pi^2 n^2}} \nonumber \eea
$$z \cot z = 1-2 \sum_{n=1}^\infty \frac{z^2}{n^2 \pi^2} \frac{1}{1-\frac{z^2}{n^2 \pi^2}}= 1-2\sum_{n=1}^\infty \frac{z^2}{n^2 \pi^2}\left(1+ \frac{z^2}{n^2 \pi^2}+ \frac{z^4}{n^4 \pi^4}+...\right)$$
\be z \cot z = 1-2\sum_{n=1}^\infty \sum_{k=1}^\infty \frac{z^{2 k}}{n^{2 k} \pi^{2 k}} \label{zeta1}\ee
From another side, we can represent this sum in terms of the Bernoulli numbers.
In the generating function (\ref{bernoullinumber})$$\frac{x}{e^x-1}= \sum_{n=0}^\infty b_n \frac{x^n}{n!},$$ where $b_{2n+1}=0 \,{\rm \, for}\,\, n\geq 1,$
by choosing $x= 2iz$ and
\bea \frac{2 i z}{e^{2 i z}-1}= \frac{z e^{-i z}}{\sin z}=\frac{z (\cos z-i \sin z)}{\sin z}= \sum_{n=0}^\infty b_n \frac{(2 i z)^n}{n!}= b_0+\sum_{k=1}^\infty b_{2 k} \frac{(2 i z)^{2 k}}{(2 k)!}.\eea
we get
\be z \cot z= 1- \sum_{k=1}^\infty b_{2 k} (-1)^{k-1} \frac{2^{2 k} z^{2 k}}{(2 k)!}.\label{zeta2}\ee
Here we used the fact that $b_{2 k+1}=0$ for $k=1,2,...$. It follows obviously from observation that l.h.s. is even function of $z.$

In this form, function on the l.h.s has infinite set of simple poles at $z= \pm \pi,\, \pm 2 \pi,...$. If $|z|< \pi,$ then it is analytic and has unique expansion to Taylor series around $z=0.$

Comparing the expressions (\ref{zeta1}) and (\ref{zeta2}),we obtain

\be \sum_{n=1}^\infty \frac{1}{n^{2 k}}= (-1)^{k-1} b_{2 k} \frac{2^{2 k-1}}{(2 k)!} \pi^{2 k}.\label{zeta3}\ee
This gives relation between different power of zeros of sine function and Bernoulli numbers.
The left-hand side of this equation is the \textbf{Riemann Zeta function}
\be \zeta(s)= \sum_{n=1}^\infty \frac{1}{n^s} \ee of even argument. This is why, we get expression of this Zeta function in terms of Bernoulli numbers
\be \zeta(2 k)= (-1)^{k-1} b_{2 k} \frac{2^{2 k-1}}{(2k)!} \pi^{2k}.\ee

The following are the first few values of the Riemann zeta function:
\be \zeta(2)=\sum_{n=1}^\infty \frac{1}{n^2}= \frac{\pi^2}{6}, \label{riemannzeta}\ee
\be \zeta(4)= \sum_{n=1}^\infty \frac{1}{n^4}=\frac{\pi^4}{90}, \label{riemannzeta2}\ee
\be \zeta(6)=\sum_{n=1}^\infty \frac{1}{n^6}= \frac{\pi^6}{945}.\ee
\section{q-Bernoulli Numbers and Zeros of q-Sine Function}
Now we are going to find the similar relation between zeros of the $q$-sine function and the $q$-Bernoulli numbers.

\subsection{$q$-Bernoulli Polynomials and Numbers}

First we introduce the $q$-analogue of Bernoulli polynomials and Bernoulli numbers. The
generating function for $q$-Bernoulli polynomials is defined in terms of Jackson's q-exponential functions as follows
\be F_x(z)= \frac{z e_q(x z)}{E_q(\frac{z}{2}) \left(e_q(\frac{z}{2})-e_q(-\frac{z}{2})\right)}= \frac{z e_q(x z) e_q(-\frac{z}{2})}{e_q(\frac{z}{2})-e_q(-\frac{z}{2})}= \sum_{n=0}^\infty B_n^q(x) \frac{z^n}{[n]!}, \ee
where the Jackson's $q$-exponential functions are \cite{Kac et al.}
\bea e_q(x)= \sum_{n=0}^\infty \frac{x^n}{[n]_q!}\,, \,\,\,\,\,\,\,\,\,\, E_q(x)= \sum_{n=0}^\infty q^{\frac{n(n-1)}{2}}\frac{x^n}{[n]_q!}, \label{qexponentials} \eea
and $[n]_q!= [1]_q [2]_q ... [n]_q, \,\,\, [n]_q= \frac{q^n-1}{q-1}.$
Two $q$-exponential functions are related to by the next formula $$e_q(x) E_q(-x)=1.$$
By $q$-differentiation the generating function with respect to $x,$ it is easy to obtain the recursion formula
\be D_q ^x B_n^q(x)=[n]_q B_{n-1}^q (x),\ee where $D_q f(x)\equiv \frac{f(q x)-f(x)}{(q-1)x}$ and $B_0^q(x)=1.$
In the limiting $q \rightarrow 1,$ this relation reduces to the standard recursion formula (\ref{differentialrelation}).

 For $n\geq 0,$ $b_n^q \equiv B_n^q(0)$ we called the \textbf{$q$-Bernoulli numbers}.

According to above definition, the generating function for $q$-Bernoulli numbers is given by
\be F_0(z)=\frac{z}{E_q(\frac{z}{2}) \left(e_q(\frac{z}{2})-e_q(-\frac{z}{2})\right)}= \sum_{n=0}^\infty b_n^q \frac{z^n}{[n]!}.\label{genfuncbernoullinum} \ee

By the definition of Jackson $q$-exponential functions (\ref{qexponentials}) we expand this generating function as
\bea & & \frac{z}{E_q(\frac{z}{2}) \left(e_q(\frac{z}{2})-e_q(-\frac{z}{2})\right)} \nonumber \\
&=& \frac{z}{\left(1+\frac{z}{2}+ q\frac{z^2}{2^2[2]!}+q^3 \frac{z^3}{2^3 [3]!}+q^6 \frac{z^4}{2^4 [4]!}+...\right) \left(z+\frac{z^3}{2^2[3]!}+\frac{z^5}{2^4 [5]!}+...\right)} \nonumber \\
&=& b_0^q+b_1^q z+ b_2^q \frac{z^2}{[2]!}+b_4^q \frac{z^4}{[4]!}+... \nonumber \\
&=& \frac{1}{1+ \frac{z}{2}+z^2\left(\frac{1}{2^2[3]!}+\frac{q}{2^2[2]!}\right) +z^3 \left(\frac{1}{2^3 [3]!}+ \frac{q^3}{2^3 [3]!}\right)+z^4\left(\frac{q}{2^4 [2]! [3]!}+ \frac{q^6}{2^4 [4]!}+ \frac{1}{2^4 [5]!}\right)+...} \nonumber \\
&=& \frac{1}{1+\frac{z}{2}+ A z^2 +B z^3 +C z^4 +...}= 1- \left(\frac{z}{2}+ A z^2 +B z^3 +C z^4 +... \right)\nonumber \\
&+& \left(\frac{z}{2}+ A z^2 +B z^3 +C z^4 +... \right)^2
-  \left(\frac{z}{2}+ A z^2 +B z^3 +C z^4 +... \right)^3 \nonumber \\
&+& \left(\frac{z}{2}+ A z^2 +B z^3 +C z^4 +... \right)^4 +...,  \eea
where $$A \equiv \frac{[4]}{2^2 [3]!},$$
$$B \equiv \frac{q^3 +1}{2^3 [3]!},$$
$$C \equiv \frac{[5] q^6+1}{2^4 [5]!}+ \frac{q}{2^4 [2]! [3]!}$$

Comparing terms with the same power of $z$ we get first few $q$-Bernoulli numbers.

For term $z^2$ we have $$-A +\frac{1}{4}= b_2^q \frac{1}{[2]!} \Rightarrow b_2^q =\frac{1}{4}\left([2]-\frac{1}{[3]}-q\right) $$
and
$$b_2^q =\frac{1}{4}\left([2]-\frac{1}{[3]}-q\right).$$
For term $z^4$ we get \bea & & -C +A^2 +B -\frac{3}{4} A + \frac{1}{16}= b_4^q \frac{1}{[4]!}\Rightarrow \nonumber \\
& & b_4^q= \frac{[4]}{2^4}\left([3]!-[2]^3+\frac{[4]^2}{[3]!}-\frac{q}{[2]!}-\frac{[5]q^6+1}{[5][4]}\right)\eea
and as a result $$ b_4^q= \frac{[4]}{2^4}\left([3]!-[2]^3+\frac{[4]^2}{[3]!}-\frac{q}{[2]!}-\frac{[5]q^6+1}{[5][4]}\right).$$

\be b_0^q=1, \,\,\,\,\, b_1^q=-\frac{1}{2},\,\,\,\,\,b_2^q= \frac{1}{4}\left([2]-\frac{1}{[3]}-q\right), \,\,\,\,\, b_3^q=0,\label{fewbernoulli} \ee
 \be b_4^q= \frac{[4]}{2^4}\left([3]!-[2]^3+\frac{[4]^2}{[3]!}-\frac{q}{[2]!}-\frac{[5]q^6+1}{[5][4]}\right). \label{qbernoulli4}\ee
By choosing $z\equiv 2it$ in generating function  (\ref{genfuncbernoullinum}), we obtain
\be F_0(2 it)= \frac{2 it}{E_q(i t)\left(e_q(i t)-e_q(-i t)\right)}= \frac{t}{E_q(it) \sin_q t}= \frac{t e_q(-i t)}{\sin_q t}.\ee
From the $q$-analogue of Euler identity $e_q(i x)=\cos_q x+i \sin_q x, $ we have
\bea F_0(2 it)= \frac{t}{ \sin_q t}(\cos_q t-i \sin_q t)= t \cot_q t-i t &=&  \sum_{n=0}^\infty b_n^q \frac{(2 i t)^n}{[n]!} \nonumber \\
&=& b_0^q + b_1^q (2 it)+  \sum_{n=2}^\infty b_n^q \frac{(2 i t)^n}{[n]!} ,\nonumber \eea
where $$\cot_q t = \frac{\cos_q t}{\sin_q t}$$ and $$\cos_q t= \frac{e_q(i t)+e_q(-i t)}{2},\,\,\,\,\,\,\,\,\,\, \sin_q t= \frac{e_q(i t)-e_q(-i t)}{2i}.$$

Then, substituting $b_0^q$ and $b_1^q$ into the above equality we get
\be t \cot_q t= 1+\sum_{n=2}^\infty b_n^q \frac{(2 i t)^n}{[n]!}, \ee or
\be t \cot_q t= 1+\sum_{k=1}^\infty b_{2 k}^q \frac{(2 i t)^{2 k}}{[2 k]!} \label{cot1}. \ee
Here the left-hand side is even function of $t,$ so that in the last sum odd coefficients vanish $b_{2 k+1}=0$ for $k=1,2,...$.\\

\subsection{Zeros of $q$-Sine Function}
Now we like to express the l.h.s. of (\ref{cot1}) in terms of zeros of $\sin_q x$ function. We start with proposition :

\begin{prop}{\textbf{$q$-Generalized Multiple Product Leibnitz Rule}:}
\bea D_q (f_1(x) f_2(x)...f_n(x))&=& \left(D_q f_1(x) \right) f_2(x)...f_n(x) \nonumber \\
&+& f_1(q x)\left(D_q f_2(x)\right) f_3(x)...f_n(x) \nonumber \\ &+& ...\nonumber \\ &+& f_1(q x) f_2(q x)...f_{n-1}(q x) \left(D_q f_n(x)\right) \eea
\end{prop}
\begin{prf}
For $n=1$ it is evident. For $n=2$ it gives $q$-Leibnitz rule \cite{Kac et al.}
$$ D_q \left(f_1(x) f_2(x)\right)= \left(D_q f_1(x)\right) f_2(x)+ f_1(q x) \left(D_q f_2(x)\right).$$
Suppose it is true for some  $n.$ Then by induction
\bea D_q (f_1(x) f_2(x)...f_n(x)f_{n+1}(x))&=&D_q (f_1(x) f_2(x)...f_n(x)) f_{n+1}(x) \nonumber \\
&+&f_1(q x) f_2(q x)...f_{n}(q x) \left(D_q f_{n+1}(x)\right) \nonumber \\ &=& \left(\left(D_q f_1(x) \right) ...f_n(x)+ ...+f_1(q x)... \left(D_q f_n(x)\right)\right) f_{n+1} (x) \nonumber \\&+& f_1(q x) f_2(q x)...f_{n}(q x) \left(D_q f_{n+1}(x)\right), \nonumber  \eea
which is the desired result.\hfill \rule{1.6ex}{1.6ex}
\end{prf}
According to the above proposition we have the following rule of differentiation multiple product \textbf{(the $q$-analogue of logarithmic derivative)}
\bea \frac{D_q(f_1 f_2...f_n)}{f_1 f_2...f_n}= \frac{f'_1(x)}{f_1(x)}+ \frac{f_1(q x)}{f_1(x)}\frac{f'_2(x)}{f_2(x)}+...+\frac{f_1(q x)}{f_1(x)}...\frac{f_{n-1}(q x)}{f_{n-1}(x)}\frac{f'_n(x)}{f_n(x)}\label{logder} \eea
\textbf{Example:} If $f_k= (x-x_k)$ and $f_1...f_n= \prod_{k=1}^n (x-x_k)$ is function with $n$ zeros, $x_1,...,x_n,$ then we have
\bea \frac{D_q \left(\prod_{k=1}^n (x-x_k)\right)}{\prod_{k=1}^n (x-x_k)}&=& \frac{1}{(x-x_1)}+ \frac{(q x-x_1)}{(x-x_1)} \frac{1}{(x-x_2)}+ \frac{(q x-x_1)}{(x-x_1)} \frac{(q x-x_2)}{(x-x_2)} \frac{1}{(x-x_3)} \nonumber \\
&+& ...+\frac{(q x-x_1)}{(x-x_1)} \frac{(q x-x_2)}{(x-x_2)} \frac{(q x-x_3)}{(x-x_3)}...\frac{(q x-x_{n-1})}{(x-x_{n-1})} \frac{1}{(x-x_n)}, \nonumber  \eea
as a simple pole expansion.\\
Expanded to simple fractions this expression can be rewritten as
\be \frac{D_q \left(\prod_{k=1}^n (x-x_k)\right)}{\prod_{k=1}^n (x-x_k)}= \sum_{k=1}^n \frac{A_k}{x-x_k},\ee where coefficients
\bea A_k &=& Res_{|{x=x_k}}\frac{D_q \left(\prod_{k=1}^n (x-x_k)\right)}{\prod_{k=1}^n (x-x_k)} \nonumber \\
&=& Res_{|{x=x_k}}\left(\frac{1}{(x-x_1)}+ \frac{(q x-x_1)}{(x-x_1)} \frac{1}{(x-x_2)}+...+ \frac{(q x-x_1)}{(x-x_1)} \frac{(q x-x_2)}{(x-x_2)} ...\frac{1}{(x-x_n)} \right) \nonumber \eea
Particularly, for $n=2,$
\bea A_1&=& \lim_{x \rightarrow x_1} \left((x-x_1) \left(\frac{1}{x-x_1}+ \frac{q x-x_1}{(x-x_1)(x-x_2)}\right) \right) \nonumber \\
&=& 1+\frac{x_1(q-1)}{(x_1-x_2)}= \frac{q x_1-x_2}{x_1-x_2},\eea
\bea A_2 &=& \lim_{x \rightarrow x_2} \left((x-x_2) \left(\frac{1}{x-x_1}+ \frac{q x-x_1}{(x-x_1)(x-x_2)}\right) \right) \nonumber \\
&=& \frac{q x_2 -x_1}{x_2-x_1},\eea
and we get
\be \frac{D_q\left((x-x_1) (x-x_2)\right)}{(x-x_1) (x-x_2)}= \left(\frac{q x_1-x_2}{x_1-x_2}\right) \frac{1}{x-x_1}+\left(\frac{q x_2-x_1}{x_2-x_1}\right) \frac{1}{x-x_2}.\ee \hfill \rule{1.6ex}{1.6ex}

We consider $\sin_q x$ function as an infinite product in terms of its zeros \\ $x_n\equiv x_n(q)$ in the following form
\be \sin_q x= x \prod _{n=1}^\infty \left(1-\frac{x^2}{x_n^2} \right)= x \left(1-\frac{x^2}{x_1^2}\right)\left(1-\frac{x^2}{x_2^2}\right) ... \label{zerosin}\ee
By using the above property (\ref{logder}), we have

\bea \frac{D_q \sin_q x}{\sin_q x}&=& \cot_q x= \frac{D_q \left(x \prod _{n=1}^\infty \left(1-\frac{x^2}{x_n^2} \right)\right)}{x \prod _{n=1}^\infty \left(1-\frac{x^2}{x_n^2} \right)} \nonumber \\ &=& \frac{1}{x}+ \frac{q x}{x} \frac{\left(-[2]\frac{x}{x_1^2}\right)}{\left(1-\frac{x^2}{x_1^2}\right)}+ \frac{q x}{x} \frac{\left(1-q^2 \frac{x^2}{x_1^2}\right)}{\left(1-\frac{x^2}{x_1^2}\right)}\frac{\left(-[2]\frac{x}{x_2^2}\right)}{\left(1-\frac{x^2}{x_2^2}\right)}+... \nonumber \\
&+& \frac{q x}{x}\frac{\left(1-q^2 \frac{x^2}{x_1^2}\right)}{\left(1-\frac{x^2}{x_1^2}\right)}\frac{\left(1-q^2 \frac{x^2}{x_2^2}\right)}{\left(1-\frac{x^2}{x_2^2}\right)}...\frac{\left(-[2]\frac{x}{x_n^2}\right)}{\left(1-\frac{x^2}{x_n^2}\right)}+... ,  \eea
where we  ordered zeros as $|x|<|x_1|<|x_2|<...<|x_n|<...,$ so that $|\frac{x}{x_k}|<1,$ for any $k.$
The above expression can be written in a compact form as follows
\be x \cot_q x= 1-[2]q \sum _{n=1}^\infty \frac{\frac{x^2}{x_n^2}}{\left(1-\frac{x^2}{x_n^2}\right)} \prod_{k=1}^{n-1} \frac{\left(1-q^2 \frac{x^2}{x_k^2}\right)}{\left(1-\frac{x^2}{x_k^2}\right)}.\label{cot2}\ee

\subsubsection{Quadratic order: }
Now we compare expressions (\ref{cot1}) and (\ref{cot2}) by equating equal powers in $x^2:$
\bea & & 1+b_2^q \frac{-4 x^2}{[2]!}+ b_4^q \frac{2^4 x^4}{[4]!}+...= \nonumber \\
& &\sum_{n=1}^\infty \frac{x^2}{x_n^2} \left(1+\frac{x^2}{x_n^2}+\left(\frac{x^2}{x_n^2}\right)^2+...\right) \cdot \left(1+(1-q^2)\frac{x^2}{x_1^2}+(1-q^2)\left(\frac{x^2}{x_1^2}\right)^2+...\right) \cdot \nonumber \\
& & \left(1+(1-q^2)\frac{x^2}{x_2^2}+(1-q^2)\left(\frac{x^2}{x_2^2}\right)^2+...\right) \cdot \nonumber \\
& & ... \nonumber \\
& &\left(1+(1-q^2)\frac{x^2}{x_{n-1}^2}+(1-q^2)\left(\frac{x^2}{x_{n-1}^2}\right)^2+...\right).\eea
At the order $x^2$  we have $$[2]q \sum_{n=1}^\infty \frac{1}{x_n^2}= b_2^q \frac{4}{[2]!}$$ and using (\ref{fewbernoulli}) for the value of Bernoulli number $b_2^q= \frac{1}{4}\left([2]-q-\frac{1}{[3]}\right),$ we obtain
\be \sum_{n=1}^\infty \frac{1}{x_n^2(q)}= \frac{1}{[3]!}.\ee
In the limiting case $q \rightarrow 1,$ $[3]!=6$ and we have $$\lim_{q \rightarrow 1} \sum_{n=1}^\infty \frac{1}{x_n^2(q)}= \frac{1}{6}.$$
Due to relation (\ref{riemannzeta}) $$\frac{1}{\pi^2} \zeta(2)= \sum_{n=1}^\infty \frac{1}{n^2 \pi^2}= \frac{1}{6},$$ it implies $$\lim_{q \rightarrow 1} x_n(q)= n \pi.$$

We found relation between zeros $x_n$ of $q$-sine function  and $b_2^q$ at order $x^2.$

\subsubsection{Quartic order: }
Now we will find relation at the order $x^4,$ this why let us call $$\frac{x^2}{x_n^2} \equiv \xi_n ,$$ then the above expression is written in terms of $\xi$ as follows
\be x \cot_q x= 1-[2]q \sum _{n=1}^\infty \frac{\xi_n}{1-\xi_n} \prod_{k=1}^{n-1} \frac{1-q^2 \xi_k}{1-\xi_k}= 1- b_2^q \frac{2^2}{[2]!}x^2+ b_4^q \frac{2^4}{[4]!}x^4+... \label{A3}\ee
For simplicity we denote $$A \equiv \sum _{n=1}^\infty \frac{\xi_n}{1-\xi_n} \prod_{k=1}^{n-1} \frac{1-q^2 \xi_k}{1-\xi_k},$$
then open form of the above expression gives
\bea A &=& \frac{\xi_1}{1-\xi_1}+\frac{\xi_2}{1-\xi_2}\frac{(1-q^2 \xi_1)}{1-\xi_1} +\frac{\xi_3}{1-\xi_3}\frac{(1-q^2 \xi_1)}{1-\xi_1}\frac{(1-q^2 \xi_2)}{1-\xi_2} \nonumber \\
&+& ...+ \frac{\xi_n}{1-\xi_n}\frac{(1-q^2 \xi_1)}{1-\xi_1}...\frac{(1-q^2 \xi_{n-1})}{1-\xi_{n-1}}+...\eea
For $|\frac{x}{x_n}|= |\xi_n|<1,$ Taylor expansion of the above expression is
\bea A&=& \xi_1(1+\xi_1+\xi_1^2+...)+\xi_2(1+\xi_2+\xi_2^2+...)\left(1+(1-q^2) \xi_1+(1-q^2) \xi_1^2+... \right)\nonumber \\
&+&\xi_3(1+\xi_3+\xi_3^2+...)\left(1+(1-q^2) \xi_1+(1-q^2) \xi_1^2+... \right) \left(1+(1-q^2) \xi_2+(1-q^2) \xi_2^2+... \right) \nonumber \\
&+& ...+ \xi_n(1+\xi_n+\xi_n^2+...)...\left(1+(1-q^2) \xi_{n-1}+(1-q^2) \xi_{n-1}^2+... \right)+...  \eea
Here we should consider just $\xi^2$ terms to collect order $x^4,$ so we denote
\bea B &=& \xi_1^2+ \xi_2^2+...+ \xi_n^2+...+ \xi_1 \xi_2 (1-q^2)+\xi_1 \xi_3 (1-q^2)+\xi_2 \xi_3 (1-q^2)+... \nonumber \\
&+& \xi_n \xi_1 (1-q^2) +\xi_n \xi_2 (1-q^2)+...+ \xi_n \xi_{n-1} (1-q^2)+... \nonumber \\
&=& \sum_{k=1}^\infty \xi_k^2 + (1-q^2) C \label{A1}, \eea
where \be C \equiv  \sum_{k=2}^\infty \xi_1 \xi_k +\sum_{k=3}^\infty \xi_2 \xi_k+...+ \sum_{k=n+1}^\infty \xi_n \xi_k+... \label{C}\ee
By $$ \sum_{k=1}^n \xi_k \equiv S_n $$ and $$\lim_{n\rightarrow \infty} S_n=S,$$ then we can write the sums as
\bea & & \sum_{k=1}^\infty \xi_k= S \nonumber \\
& &\sum_{k=2}^\infty \xi_k = \sum_{k=1}^\infty \xi_k -\xi_1 = S-S_1 \nonumber \\
& &\sum_{k=3}^\infty \xi_k = \sum_{k=1}^\infty \xi_k -\xi_1-\xi_2 = S-S_2 \nonumber \\
& &... \nonumber \\
& &\sum_{k=n}^\infty \xi_k = S-S_{n-1} \eea
Rewriting (\ref{A1}) in terms of $S,$ we obtain
\bea B &=& S^2 + (1-q^2) \left(\xi_1(S-S_1)+ \xi_2 (S-S_2)+...+ \xi_n(S-S_n)+... \right) \nonumber \\
&=& S^2 + (1-q^2)\left( S^2- \xi_1 \xi_1-\xi_2(\xi_1+\xi_2)-...- \xi_n(\xi_1+\xi_2+...+\xi_n) \right)  \nonumber \\
&=& S^2 + (1-q^2)\left( S^2- \sum_{k=1}^\infty \xi_k^2 +D \right), \label{A2}\eea
where
\be D \equiv -\xi_2 \xi_1 -\xi_3 (\xi_1+\xi_2)-...-\xi_n(\xi_1+\xi_2+...+\xi_{n-1})+...\ee
or
\bea D&=& -\xi_1(\xi_2+\xi_3+...+ \xi_n)-\xi_2(\xi_3+\xi_4+...+\xi_n)-...-\xi_n(\xi_{n+1}+...)-... \nonumber \\
&=& -\xi_1(S-S_1)-\xi_2(S-S_2)-...-\xi_n(S-S_n).\eea
Comparing with (\ref{C}) we find $D=-C,$ then  by equating (\ref{A1}) and (\ref{A2})
\be B=S^2 + (1-q^2)\left( S^2- \sum_{k=1}^\infty \xi_k^2 -C  \right)= S^2 +(1-q^2)C  ,\ee
we get \be C= \frac{1}{2} S^2 -\frac{1}{2} \sum_{k=1}^\infty \xi_k^2 .\ee
It gives
\bea B &=& \sum_{k=1}^\infty \xi_k^2 + (1-q^2)C \nonumber \\
&=&\sum_{k=1}^\infty \xi_k^2 + (1-q^2)\left(\frac{1}{2} S^2 -\frac{1}{2} \sum_{k=1}^\infty \xi_k^2\right) \nonumber \\
&=& \left(1+\frac{q^2-1}{2}\right) \sum_{k=1}^\infty \xi_k^2 -\left(\frac{q^2-1}{2}\right) S^2.\eea
For $x^4$ term then we have
\be b_4^q \frac{2^4}{[4]!}x^4= -[2]q B,\ee
and substituting $\xi_k= \frac{x^2}{x_k^2}$ and $S= \sum_{k=1}^\infty \xi_k = \sum_{k=1}^\infty \frac{x^2}{x_k^2}$ in $B,$
finally we obtain
\be [2]q \left( 1+\frac{q^2-1}{2}\right) \sum_{k=1}^\infty \frac{1}{x_k^4}= \frac{8(q^2-1)}{[2]^3 q} (b_2^q)^2- \frac{16}{[4]!}b_4^q,\ee
where
$$b_2^q= \frac{1}{4}\left([2]-\frac{1}{[3]}-q\right),$$
$$b_4^q= \frac{[4]}{2^4}\left([3]!-[2]^3+\frac{[4]^2}{[3]!}-\frac{q}{[2]!}-\frac{[5]q^6+1}{[5][4]}\right). $$

In the limiting case $q \rightarrow 1,$
$$\lim_{q \rightarrow 1} \sum_{n=1}^\infty \frac{1}{x_n^4(q)}= \frac{1}{90}.$$
From the relation (\ref{riemannzeta2}) we get
$$\frac{1}{\pi^4} \zeta(4)= \sum_{n=1}^\infty \frac{1}{n^4 \pi^4}= \frac{1}{90}.$$

This procedure can be continued to higher order relating higher powers of zeros of q-sine function with q-Bernoulli numbers. However, procedure at every next step become more and more complicated. Moreover, it is another problem finding exact form of zeros $x_n(q)$ of $\sin_q x$ from infinite set of equations.

\section{Conclusions}

In this paper, we have introduced the set of $q$-Bernoulli numbers and found relations between these numbers and zeros of $q$-sine function. The exact formula for these zeros is not known, this is why, our relations can be useful in analyzing asymptotic formulas and approximate form of the zeros.

Here we mention some simple approximate formula for zeros of $q$-sine function. It is coming from analyzing graph of $\sin_q x$ and for $q>1$ is
$$x_2= q^2 x_1,\,\,\,\,
x_3= q^2 x_2= q^4 x_1,...,x_n= q^{2(n-1)} x_1,$$
then,
\bea \frac{1}{[3]!}&=& \sum_{n=1}^\infty \frac{1}{x_n^2}= \frac{1}{x_1^2} + \frac{1}{x_2^2}+...= \frac{1}{x_1^2}+ \frac{1}{q^4 x_1^2}+ \frac{1}{q^8 x_1^2}+... \nonumber \\
&=& \frac{1}{x_1^2} \left(1+ \frac{1}{q^4}+ \left(\frac{1}{q^4}\right)^2+... \right)\eea
and from sum of  geometric series in $\frac{1}{q^4}$ we have $$\frac{1}{[2] [3]}= \frac{1}{x_1^2} \frac{q^4}{q^4-1}.$$
From this expression we have the first root as $$x_1= \pm \sqrt{[2][3]\frac{q^4}{q^4-1}}.$$

As a result, (\ref{zerosin}) can be written in the following form
\be \sin_q x= x \prod_{n=1}^\infty \left( 1- \frac{[4](q-1) x^2}{q^{4 n} [3]!} \right),\ee
where for wave number we have the discrete set $x_n^2 = \frac{q^{4 n}[2][3]}{[4](q-1)}.$

These results can be useful in solving BVP for $q$-wave equation and $q$-Sturm-Liouville problem. These questions are under investigation now.

 \section*{Acknowledgments}
This work was support by
TUBITAK (The Scientific and Technological Research Council of Turkey), TBAG Project 110T679 and Izmir Institute of
Technology. And one of the authors (S.Nalci) was partially supported by TUBITAK scholarship of graduate students.


\begin{thebibliography}{xx}


\bibitem{Kac et al.} V. Kac and P. Cheung,\newblock Quantum Calculus, Springer, New York, 2002.
\bibitem{Exton} H. Exton, \newblock q-Hypergeometric Functions and Applications, John Wiley and Sons, 1983.
\bibitem{Jackson} F.H. Jackson ,\newblock A Basic Sine and Cosine with Symbolic Solutions of certain Differential Equations, Proc. Edin. Math. Soc. 22, 28-39, 1904.
\bibitem{Sury} B. Sury ,\newblock Bernoulli Numbers and the Riemann Zeta Function, Resonance,  2003.
\bibitem{Al-Salam} W. A. Al-Salam,\newblock q-Bernoulli numbers and polynomials, Math. Nachr. Vol. 17, pp. 239-260,  1959.
\bibitem{Carlitz} L. Carlitz,\newblock q-Bernoulli numbers and polynomials, Duke Math. Jour., Vol. 15,  pp. 987-1000,  1948.
\bibitem{Ryoo} C.S. Ryoo, \newblock A note on q-Bernoulli numbers and polynomials, Applied Mathematics Letters Vol. 20,  pp. 524-531,  2007.
\end{thebibliography}
\end{document}